\documentclass[11pt,reqno]{amsproc}
\usepackage{amsmath, amsthm, amsfonts, amssymb}
\usepackage{enumitem}
\usepackage[bookmarksnumbered, colorlinks, plainpages]{hyperref}
\hypersetup{pdftitle={Relative nuclear dimension of discrete groups: exactness, K-theory and holonomy},pdfauthor={Ali Jabbari},pdfsubject={Operator algebras and geometric group theory}}
\newcommand{\dimnuc}{\operatorname{dim}_{\mathrm{nuc}}}
\newcommand{\rdim}{\operatorname{dim}_{\mathrm{nuc}}^{\mathrm{u}}}
\newcommand{\Cr}{C^*_r}
\newcommand{\Cu}{C^*_u}
\newcommand{\asdim}{\operatorname{asdim}}

\newtheorem{theorem}{\textbf{Theorem}}[section]
\newtheorem{corollary}[theorem]{\bf{Corollary}}
\newtheorem{lemma}[theorem]{\bf{Lemma}}
\newtheorem{proposition}[theorem]{\bf{Proposition}}
\newtheorem{question}[theorem]{\bf{Question}}
\newtheorem*{maintheorem}{\bf{Main Theorem}}
\theoremstyle{definition}
\newtheorem{definition}[theorem]{\bf{Definition}}
\newtheorem{remark}[theorem]{\bf{Remark}}
\newtheorem{example}[theorem]{\bf{Example}}
\numberwithin{equation}{section}

\begin{document}
\sloppy

\title[Relative nuclear dimension of discrete groups]{Relative nuclear dimension of discrete groups: exactness, $K$-theory and holonomy}

\author{Ali Jabbari}
\address{Kyrgyz-Turkish Manas University, Bishkek, Kyrgyzstan}
\email{ali.jabbari@manas.edu.kg}

\begin{abstract}
For a countable discrete group $\Gamma$, let $\rdim(\Gamma)$ be the nuclear dimension of the canonical inclusion $\Cr(\Gamma)\hookrightarrow\Cu(\Gamma)$. We show that it takes values in $\{0,1,\infty\}$, is finite exactly for exact groups, and can vanish only for amenable groups. We give two proofs of the upper bound for exact groups, one of which uses Cantor crossed products and is independent of the recent theorem of Guo, Li, Zhang and Zhu \cite{GLZZ2026} on uniform Roe algebras. A $K$-theoretic obstruction shows that $\rdim(\Gamma)=1$ whenever $\Gamma$ is amenable and its abelianization is not torsion. In particular, we compute the invariant for all countable abelian groups. We show that the regular inclusion of $\mathbb Z^n$ is injective on $K$-theory, and deduce that for a group with a finite-index normal subgroup $\mathbb Z^n$ and holonomy group $F$ the image of $K_1(\iota_\Gamma)$ is free abelian of rank $\sum_k\dim_{\mathbb Q}(\Lambda^{2k+1}\mathbb Q^n)^F$. So $\rdim(\Gamma)=1$ whenever some odd exterior power of the holonomy representation has nonzero invariants. This produces virtually abelian groups with finite abelianization and relative nuclear dimension one, among them an index-two subgroup of $D_\infty^3$, while every nontrivial finitely generated torsion-free virtually abelian group has relative nuclear dimension one. In the opposite direction, an AF crossed-product theorem of Bratteli, Evans and Kishimoto yields $\rdim(D_\infty)=0$, although $\dimnuc(\Cu(D_\infty))=1$. Thus the relative invariant is neither monotone under passage to subgroups nor invariant under finite-index inclusions or coarse equivalence. We also prove permanence of dimension zero under finite direct products and increasing unions, obtaining examples of arbitrarily large finite asymptotic dimension for which the regular inclusion has dimension zero.
\end{abstract}

\maketitle

\section{Introduction}

The nuclear dimension of Winter and Zacharias \cite{WZ2010} is defined through approximations of the identity map by completely positive maps through finite-dimensional algebras, so it is infinite for every non-nuclear C$^*$-algebra. Reduced group C$^*$-algebras of non-amenable groups, such as $\Cr(\mathbb F_2)$, are never nuclear. Weakening the topology of the approximation does not help: if the identity of $A$ is a point--weak$^*$ limit in $A^{**}$ of maps through finite-dimensional algebras whose back maps are uniformly bounded sums of a fixed finite number of order zero maps and whose incoming maps are uniformly completely bounded, then $A$ is nuclear (Proposition~\ref{prop:obstruction}). So one should keep norm approximations and change the map being approximated.

For a group $\Gamma$ there is a canonical C$^*$-algebra containing $\Cr(\Gamma)$ which is nuclear whenever $\Gamma$ is exact: the uniform Roe algebra $\Cu(\Gamma)=\ell^\infty(\Gamma)\rtimes_r\Gamma$ \cite{O2000, BOB}. We define the \emph{relative nuclear dimension} of $\Gamma$ as the nuclear dimension of the inclusion map,
\[
\rdim(\Gamma)=\dimnuc\big(\iota_\Gamma\colon\Cr(\Gamma)\hookrightarrow\Cu(\Gamma)\big).
\]

Trivially $\rdim(\Gamma)\le\dimnuc(\Cu(\Gamma))$. Winter and Zacharias proved $\dimnuc(\Cu(X))\le\asdim X$ and asked whether equality holds \cite[Question~9.5]{WZ2010}. Guo, Li, Zhang and Zhu have answered this in the negative: $\dimnuc(\Cu(X))\le1$ for every space $X$ of bounded geometry with Property A, so $\dimnuc(\Cu(X))\in\{0,1,\infty\}$, with $0$ exactly when $\asdim X=0$ by Li and Willett \cite{GLZZ2026, LW2018}. For groups this reads
\[
\dimnuc(\Cu(\Gamma))=\begin{cases}0,&\Gamma\ \text{locally finite},\\ 1,&\Gamma\ \text{exact, not locally finite},\\ \infty,&\Gamma\ \text{not exact}.\end{cases}
\]

The question is how much of this dimension is already needed for the subalgebra $\Cr(\Gamma)$.

\begin{maintheorem}
Let $\Gamma$ be a countable discrete group.
\begin{enumerate}[label=(\roman*)]
\item $\rdim(\Gamma)\in\{0,1,\infty\}$, and it is finite if and only if $\Gamma$ is exact.
\item If $\rdim(\Gamma)=0$, then $\Gamma$ is amenable and $K_1(\iota_\Gamma)=0$, and the abelianization $\Gamma_{\mathrm{ab}}$ is a torsion group.
\item If $\Gamma$ is locally finite, then $\rdim(\Gamma)=0$. If $\Gamma$ is exact and non-amenable, or amenable with non-torsion abelianization, then $\rdim(\Gamma)=1$. The same holds if $\Gamma$ has a normal subgroup of finite index isomorphic to $\mathbb Z^n$ whose holonomy group $F\subset\mathrm{GL}(n,\mathbb Z)$ satisfies $(\Lambda^{2k+1}\mathbb Q^n)^F\ne0$ for some $k\ge0$.
\item For countable abelian $\Gamma$, $\rdim(\Gamma)=0$ if $\Gamma$ is torsion and $\rdim(\Gamma)=1$ otherwise.
\item For the infinite dihedral group $D_\infty=\mathbb Z\rtimes_{-1}\mathbb Z/2\mathbb Z$, one has
\[
\rdim(D_\infty)=0,\qquad \rdim(\mathbb Z)=1,\qquad \dimnuc(\Cu(D_\infty))=1.
\]
So $\rdim$ is not subgroup-monotone, is not invariant under finite-index passage, and is not a coarse invariant.
\item The class of countable groups of relative nuclear dimension zero is closed under finite direct products and increasing unions. In particular, $\rdim(D_\infty^n)=0$ and $\asdim(D_\infty^n)=n$ for each $n\ge1$.
\item The index-two subgroup $\Gamma_0=\mathbb Z^3\rtimes\{\mathrm{diag}(\epsilon_1,\epsilon_2,\epsilon_3):\epsilon_i=\pm1,\ \epsilon_1\epsilon_2\epsilon_3=1\}$ of $D_\infty^3$ has $(\Gamma_0)_{\mathrm{ab}}\cong(\mathbb Z/2\mathbb Z)^5$ and $\rdim(\Gamma_0)=1$. So amenability together with torsion abelianization does not imply dimension zero. The same holds for $\Gamma_0\times D_\infty$, which contains $\mathbb Z^4$ with finite index and has abelianization $(\mathbb Z/2\mathbb Z)^7$.
\item The inclusion $\Cr(\mathbb Z^n)\hookrightarrow\Cu(\mathbb Z^n)$ is injective on $K_0$ and $K_1$. If $\Gamma$ has a normal subgroup $H\cong\mathbb Z^n$ of finite index with holonomy group $F$, then the image of $K_1(\iota_\Gamma)$ is free abelian of rank $\sum_{k\ge0}\dim_{\mathbb Q}(\Lambda^{2k+1}\mathbb Q^n)^F$. In particular $K_1(\iota_\Gamma)=0$ if and only if no odd exterior power of $\mathbb Q^n$ has nonzero $F$-invariants, and $\rdim(\Gamma)=1$ for every nontrivial finitely generated torsion-free virtually abelian group $\Gamma$.
\end{enumerate}
\end{maintheorem}

The general theory of zero-dimensional morphisms has been developed by Castillejos and Neagu \cite{CN2024}. We recall the finite-dimensional approximation mechanism and its consequences for quasidiagonality and $K$-theory, giving self-contained proofs of the facts needed here. Our emphasis is the application to the regular inclusion. The upper bound in (i) follows from \cite{GLZZ2026}; we also give an independent argument using Cantor crossed products \cite{RS2012,CJKMST2018}. The lower bound in (iii) for groups with non-torsion abelianization uses the Pimsner--Voiculescu boundary map at finitely generated stages of the crossed product. The opposing dihedral example relies on the AF theorem of Bratteli, Evans and Kishimoto \cite{BEK1993} and the orbit-map factorization of the regular inclusion. For virtually abelian groups the $K_1$-obstruction can be computed completely. The regular inclusion of $\mathbb Z^n$ is injective on $K$-theory, detected by iterated Pimsner--Voiculescu boundary maps and an invariant mean, and a transfer argument identifies the image of $K_1(\iota_\Gamma)$ rationally with the invariants of the holonomy group in the odd exterior powers. This shows that torsion abelianization does not suffice for dimension zero: the orientation-preserving index-two subgroup of $D_\infty^3$ has finite abelianization and relative nuclear dimension one. So neither amenability nor local finiteness characterizes dimension zero, the presence of an infinite cyclic subgroup does not force dimension one, and torsion abelianization does not force dimension zero. For virtually abelian groups, the remaining question is whether the vanishing of these holonomy invariants already forces dimension zero (Question~\ref{q:virtuallyabelian}).

\section{Preliminaries}

\subsection{Nuclear dimension of maps}
A c.p.\ map $\psi\colon E\to B$ has \emph{order zero} if $\psi(x)\psi(y)=0$ whenever $x,y\in E_+$ and $xy=0$. If $E=E^{(0)}\oplus\cdots\oplus E^{(n)}$ is finite-dimensional, a c.p.\ map $\psi\colon E\to B$ is \emph{$n$-decomposable} if each restriction $\psi|_{E^{(i)}}$ is c.p.c.\ and of order zero. By the structure theorem \cite{WZ2009}, a c.p.c.\ order zero map $\psi\colon E\to B$ has the form $\psi(x)=h\,\pi(x)$, where $C=C^*(\psi(E))$, $h=\psi(1_E)\in C_+$, and $\pi\colon E\to M(C)\cap\{h\}'$ is a $*$-homomorphism. The composition of an $n$-decomposable map with a $*$-homomorphism is $n$-decomposable.

\begin{definition}[\cite{WZ2010, TW2014}]
Let $\theta\colon A\to B$ be a c.p.c.\ map. Then $\dimnuc(\theta)\le n$ if for every finite set $\mathcal F\subset A$ and every $\varepsilon>0$ there are a finite-dimensional C$^*$-algebra $E$, a c.p.c.\ map $\varphi\colon A\to E$ and an $n$-decomposable map $\psi\colon E\to B$ with $\|\psi\varphi(a)-\theta(a)\|<\varepsilon$ for $a\in\mathcal F$. The nuclear dimension of $A$ is $\dimnuc(A)=\dimnuc(\mathrm{id}_A)$.
\end{definition}

\begin{lemma}\label{lem:factor}
Let $\iota\colon A\to C$ be a c.p.c.\ map that factors as $\iota=\Theta\circ j$, where $j\colon A\to B$ is c.p.c.\ and $\Theta\colon B\to C$ is a $*$-homomorphism. Then $\dimnuc(\iota)\le\dimnuc(j)\le\dimnuc(B)$.
\end{lemma}

\begin{proof}
Given an approximation $(E,\varphi,\psi)$ for $j$, the triple $(E,\varphi,\Theta\psi)$ is one for $\iota$, since $\Theta$ preserves $n$-decomposability and is contractive. Given an approximation $(E,\varphi,\psi)$ for $\mathrm{id}_B$ on the finite set $j(\mathcal F)$, the triple $(E,\varphi j,\psi)$ is one for $j$.
\end{proof}

\subsection{The obstruction}
The next proposition explains why we approximate an inclusion rather than the identity. It is the argument behind the classical fact that the canonical inclusion $A\hookrightarrow A^{**}$ is weakly nuclear only for nuclear $A$.

\begin{proposition}\label{prop:obstruction}
Let $A$ be a C$^*$-algebra. Suppose there are constants $C,N$ and a net of finite-dimensional C$^*$-algebras $E_\lambda$ with linear maps $\varphi_\lambda\colon A\to E_\lambda$, $\|\varphi_\lambda\|_{\mathrm{cb}}\le C$, and $\psi_\lambda\colon E_\lambda\to A^{**}$ such that:
\begin{enumerate}[label=(\alph*)]
\item each $\psi_\lambda$ is a sum of at most $N$ maps of the form $x\mapsto v\,k\,\sigma(x)\,k$, with $\sigma\colon E_\lambda\to A^{**}$ a $*$-homomorphism, $k\in A^{**}_+\cap\sigma(E_\lambda)'$, $v\in A^{**}$ a partial isometry, and $\|k\|^2\le C$;
\item $\psi_\lambda\varphi_\lambda(a)\to a$ weak$^*$ in $A^{**}$ for every $a\in A$.
\end{enumerate}
Then $A$ is nuclear.
\end{proposition}

\begin{proof}
Let $B$ be a C$^*$-algebra, $x\in A\odot B$, and let $\pi\colon A\to B(H)$, $\sigma_B\colon B\to B(H)$ be representations with commuting ranges. Let $\pi^{**}$ be the normal extension of $\pi$ to $A^{**}$; its range commutes with $\sigma_B(B)$. For each summand $x\mapsto vk\sigma(x)k$ of $\psi_\lambda$, the map $\rho_\lambda(e\otimes b)=\pi^{**}(\sigma(e))\sigma_B(b)$ extends to a $*$-homomorphism on $E_\lambda\otimes B$ (whose C$^*$-norm is unique). The map $y\mapsto\pi^{**}(v)\pi^{**}(k)\rho_\lambda(y)\pi^{**}(k)$ is therefore bounded by $\|k\|^2\le C$. The map $\varphi_\lambda\otimes\mathrm{id}_B\colon A\otimes_{\min}B\to E_\lambda\otimes B$ has norm at most $C$. Composing, we obtain maps $\Theta_\lambda\colon A\otimes_{\min}B\to B(H)$ with $\|\Theta_\lambda\|\le NC^2$ and $\Theta_\lambda(a\otimes b)=\pi^{**}(\psi_\lambda\varphi_\lambda(a))\sigma_B(b)$. For vectors $\xi,\eta$ the functional $y\mapsto\langle\pi^{**}(y)\sigma_B(b)\xi,\eta\rangle$ is normal on $A^{**}$, so $\langle(\pi\times\sigma_B)(x)\xi,\eta\rangle=\lim_\lambda\langle\Theta_\lambda(x)\xi,\eta\rangle$, and $\|x\|_{\max}\le NC^2\|x\|_{\min}$. The canonical surjection $A\otimes_{\max}B\to A\otimes_{\min}B$ is then injective, hence isometric, for every $B$, and $A$ is nuclear \cite{BOB}.
\end{proof}

Every c.p.c.\ order zero map from a finite-dimensional algebra has the form in (a) with $v=1$ and $k=h^{1/2}$, by the structure theorem. So approximations of the identity in the style of nuclear dimension cannot be finite for non-nuclear algebras, even after weakening the topology to the weak$^*$ topology of the bidual and the downward maps to completely bounded ones.

\subsection{The uniform Roe algebra}
Let $\Gamma$ be a countable discrete group with left regular representation $\lambda$ on $\ell^2(\Gamma)$. Let $\Cu(\Gamma)$ be the C$^*$-algebra on $\ell^2(\Gamma)$ generated by the multiplication operators $\ell^\infty(\Gamma)$ and by $\lambda(\Gamma)$. It is the uniform Roe algebra of $\Gamma$ with a proper right-invariant metric, it does not depend on the choice of that metric, and it is canonically isomorphic to the reduced crossed product $\ell^\infty(\Gamma)\rtimes_r\Gamma$ for the left translation action, with $u_g\mapsto\lambda_g$ \cite[Section~5.1]{BOB}. We write $\iota_\Gamma\colon\Cr(\Gamma)\hookrightarrow\Cu(\Gamma)$ for the inclusion.

We use the following facts.
\begin{enumerate}[label=(F\arabic*)]
\item $\Gamma$ is exact if and only if $\Cr(\Gamma)$ is exact, if and only if $\Cu(\Gamma)$ is nuclear \cite{KW1999, O2000, BOB}.
\item A C$^*$-algebra is exact if it admits a faithful nuclear $*$-homomorphism into $B(H)$ \cite{K1995, BOB}.
\item For a countable discrete metric space $X$ of bounded geometry, $\dimnuc(\Cu(X))\le1$ if $X$ has Property A \cite{GLZZ2026}, and $\dimnuc(\Cu(X))=0$ if and only if $\asdim X=0$ \cite{LW2018}. Earlier, Winter and Zacharias proved $\dimnuc(\Cu(X))\le\asdim X$ \cite{WZ2010}. A countable group has Property A if and only if it is exact \cite{HR2000, O2000}, and it has asymptotic dimension $0$ if and only if it is locally finite.
\item A separable unital C$^*$-algebra $A$ is quasidiagonal if for every finite $\mathcal F\subset A$ and $\varepsilon>0$ there is a c.p.c.\ map $\varphi\colon A\to M_k$ with $\|\varphi(ab)-\varphi(a)\varphi(b)\|<\varepsilon$ and $\|\varphi(a)\|>\|a\|-\varepsilon$ for $a,b\in\mathcal F$ \cite{V1991, BOB}.
\item (Rosenberg) If $\Cr(\Gamma)$ is quasidiagonal, then $\Gamma$ is amenable \cite{H1987}.
\end{enumerate}

\section{Definition and upper bounds}

\begin{definition}
The \emph{relative nuclear dimension} of a countable discrete group $\Gamma$ is
\[
\rdim(\Gamma)=\dimnuc(\iota_\Gamma),\qquad \iota_\Gamma\colon\Cr(\Gamma)\hookrightarrow\Cu(\Gamma).
\]
\end{definition}

\begin{proposition}\label{prop:crossed}
Let $\Gamma$ act on a nonempty compact Hausdorff space $X$. Then $\rdim(\Gamma)\le\dimnuc(C(X)\rtimes_r\Gamma)$.
\end{proposition}

\begin{proof}
Fix $x_0\in X$ and define $\Phi\colon C(X)\to\ell^\infty(\Gamma)$ by $\Phi(f)(g)=f(gx_0)$. With the actions $(g\cdot f)(x)=f(g^{-1}x)$ and $(g\cdot\xi)(h)=\xi(g^{-1}h)$ we have $\Phi(g\cdot f)(h)=f(g^{-1}hx_0)=(g\cdot\Phi(f))(h)$, so $\Phi$ is a unital equivariant $*$-homomorphism. It induces a $*$-homomorphism $\Theta\colon C(X)\rtimes_r\Gamma\to\ell^\infty(\Gamma)\rtimes_r\Gamma\cong\Cu(\Gamma)$ with $\Theta(u_g)=\lambda_g$ \cite[Section~4.1]{BOB}. The map $j\colon\Cr(\Gamma)\to C(X)\rtimes_r\Gamma$, $\lambda_g\mapsto u_g$, is an injective $*$-homomorphism and $\Theta\circ j=\iota_\Gamma$. Apply Lemma~\ref{lem:factor}.
\end{proof}

\begin{corollary}\label{cor:upper}
\begin{enumerate}[label=(\alph*)]
\item $\rdim(\Gamma)\le\dimnuc(\Cu(\Gamma))$.
\item $\rdim(\Gamma)\le\dimnuc(\Cr(\Gamma))$ for every countable discrete group $\Gamma$.
\item If $\Gamma$ is locally finite, then $\rdim(\Gamma)=0$.
\end{enumerate}
\end{corollary}

\begin{proof}
(a) Take $X=\beta\Gamma$ in Proposition~\ref{prop:crossed}, so that $C(X)\rtimes_r\Gamma=\ell^\infty(\Gamma)\rtimes_r\Gamma=\Cu(\Gamma)$; or apply Lemma~\ref{lem:factor} to $\iota_\Gamma=\mathrm{id}\circ\iota_\Gamma$. (b) Take $X$ a point. (c) $\Cr(\Gamma)$ is then an AF algebra, which has nuclear dimension zero \cite{WZ2010}; apply (b).
\end{proof}

\begin{lemma}\label{lem:normalize}
Let $\theta\colon A\to B$ be a unital c.p.\ map between unital C$^*$-algebras with $\dimnuc(\theta)<\infty$. Then $\theta$ is nuclear: it is a point-norm limit of maps $\psi'\varphi'$ with $\varphi'\colon A\to E'$ and $\psi'\colon E'\to B$ c.p.c.\ and $E'$ finite-dimensional.
\end{lemma}

\begin{proof}
Let $n=\dimnuc(\theta)$, $\mathcal F\subset A$ finite with $1\in\mathcal F$, and $\delta>0$. Take $(E,\varphi,\psi)$ with $\|\psi\varphi(a)-\theta(a)\|<\delta$ on $\mathcal F$; then $\|\psi\varphi(1)-1\|<\delta$. Let $e=\varphi(1)\in E_+$, $p$ its support projection, and $e^{-1/2}$ the inverse of $e^{1/2}$ in $pEp$. Put $E'=pEp$, $\varphi'(a)=e^{-1/2}\varphi(a)e^{-1/2}$ and $\psi_0(x)=\psi(e^{1/2}xe^{1/2})$. For $a\ge0$ we have $0\le\varphi(a)\le\|a\|e$, so $(1-p)\varphi(a)(1-p)=0$, hence $\varphi(a)^{1/2}(1-p)=0$ and $\varphi(a)=p\varphi(a)p$; by linearity $\varphi(A)\subseteq pEp$. So $\varphi'$ is a unital c.p.\ map into $E'$ and $\psi_0\varphi'=\psi\varphi$. The map $\psi_0$ is c.p.\ with $\|\psi_0\|=\|\psi_0(1)\|=\|\psi\varphi(1)\|<1+\delta$. With $\psi'=(1+\delta)^{-1}\psi_0$ we get $\|\psi'\varphi'(a)-\theta(a)\|<\delta+\delta\|a\|$ for $a\in\mathcal F$.
\end{proof}

\begin{proposition}\label{prop:exact}
$\rdim(\Gamma)<\infty$ implies that $\Gamma$ is exact.
\end{proposition}

\begin{proof}
By Lemma~\ref{lem:normalize}, $\iota_\Gamma$ is nuclear, and so is its composition with the inclusion $\Cu(\Gamma)\subset B(\ell^2(\Gamma))$. This is a faithful nuclear representation of $\Cr(\Gamma)$, so $\Cr(\Gamma)$ is exact by (F2), and $\Gamma$ is exact by (F1).
\end{proof}

\section[Dimension zero: quasidiagonality and K-theory]{Dimension zero: quasidiagonality and $K$-theory}

We first isolate the local finite-dimensional characterization of a zero-dimensional unital morphism (see also \cite{CN2024}). The characterization does not require separability of the codomain.

\begin{lemma}\label{lem:localfd}
Let $\theta:A\to B$ be a unital $*$-homomorphism of unital C$^*$-algebras. The following are equivalent:
\begin{enumerate}[label=(\roman*)]
\item $\dimnuc(\theta)=0$.
\item For each finite $\mathcal F\subset A$ and $\varepsilon>0$, there is a unital finite-dimensional C$^*$-subalgebra $F\subset B$ such that $\operatorname{dist}(\theta(a),F)<\varepsilon$ for every $a\in\mathcal F$.
\end{enumerate}
\end{lemma}
\begin{proof}
For (i) implies (ii), include $1_A$ in the finite set and choose a c.p.c.\ approximation $\psi\varphi$ within $\delta<1/2$ of $\theta$, with $\psi:E\to B$ order zero. Set $h=\psi(1_E)$. As $\varphi(1_A)\le1_E$, we have $(1-\delta)1_B\le h\le1_B$. By the order zero structure theorem, $\psi(x)=h\pi(x)$, where the supporting homomorphism $\pi$ takes values in $B$, since $h$ is invertible and $1_B$ belongs to $C^*(h)\subset C^*(\psi(E))$. It is unital, because $h(1_B-\pi(1_E))=\psi(1_E)-\psi(1_E)=0$ and $h$ is invertible. Thus $\|\pi(x)-\psi(x)\|\le\delta\|x\|$. The unital algebra $F=\pi(E)$ satisfies
\[
\operatorname{dist}(\theta(a),F)\le\|\theta(a)-\pi\varphi(a)\|<\delta(1+\|a\|).
\]

Choose $\delta$ sufficiently small. Conversely, for a unital finite-dimensional $F\subset B$ there is a u.c.p.\ conditional expectation $P:B\to F$, by injectivity of finite-dimensional C$^*$-algebras (or the Arveson extension theorem, applied to the matrix summands). With $j_F:F\hookrightarrow B$, we have
\[
\|j_FP\theta(a)-\theta(a)\|\le2\operatorname{dist}(\theta(a),F).
\]

The map $j_F$ is a $*$-homomorphism, hence c.p.c.\ order zero, and the desired dimension-zero approximation follows.
\end{proof}

\begin{proposition}\label{prop:qd}
Let $A\subset B$ be a unital inclusion of unital C$^*$-algebras with $A$ separable. If the inclusion map has nuclear dimension zero, then $A$ is quasidiagonal.
\end{proposition}

\begin{proof}
Let $\mathcal F\subset A$ be finite with $1\in\mathcal F$, put $\mathcal G=\mathcal F\cup\{ab:a,b\in\mathcal F\}$ and let $0<\varepsilon<1/2$. Take a c.p.c.\ map $\varphi\colon A\to E$ and a c.p.c.\ order zero map $\psi\colon E\to B$ with $\|\psi\varphi(a)-a\|<\varepsilon$ for $a\in\mathcal G$. Write $\psi(x)=h\pi(x)$ as in Section~2, with $h=\psi(1_E)$ and $C=C^*(\psi(E))$.

Since $\varphi(1)\le1_E$, we have $\psi\varphi(1)\le h\le1$, and $\|\psi\varphi(1)-1\|<\varepsilon$ gives $h\ge(1-\varepsilon)1$. So $\sigma(h)\subset[1-\varepsilon,1]$ and $\|h-1\|\le\varepsilon$. Continuous functional calculus with a function equal to $1$ on $\sigma(h)$ and $0$ at $0$ shows $1_B\in C^*(h)\subset C$. So $C$ is unital, $M(C)=C\subset B$, and $\pi$ is a $*$-homomorphism $E\to B$ commuting with $h$. Then
\[
\|\pi(x)-\psi(x)\|=\|(1-h)\pi(x)\|\le\varepsilon\|x\|,\qquad \|\pi\varphi(a)-a\|<\varepsilon(1+\|a\|)\quad(a\in\mathcal G).
\]

Let $z$ be the central projection of $E$ with $\ker\pi=(1-z)E$, put $E'=zE$ and $\varphi'=z\varphi(\cdot)$. Then $\pi|_{E'}$ is injective, hence isometric, and $\pi\varphi=\pi\varphi'$. For $a,b\in\mathcal F$,
\[
\|\varphi'(ab)-\varphi'(a)\varphi'(b)\|=\|\pi\varphi(ab)-\pi\varphi(a)\pi\varphi(b)\|\le\|\pi\varphi(ab)-ab\|+\|ab-\pi\varphi(a)\pi\varphi(b)\|,
\]
which is at most $c_{\mathcal F}\,\varepsilon$ for a constant $c_{\mathcal F}$ depending only on $\max_{a\in\mathcal F}\|a\|$. Also $\|\varphi'(a)\|=\|\pi\varphi(a)\|>\|a\|-\varepsilon(1+\|a\|)$. Embedding $E'$ unitally into some $M_k$, (F4) shows that $A$ is quasidiagonal.
\end{proof}

\begin{proposition}[$K_1$ obstruction]\label{prop:kone}
Let $A$ and $B$ be unital C$^*$-algebras and let $\theta:A\to B$ be a unital $*$-homomorphism. If $\dimnuc(\theta)=0$, then $K_1(\theta)=0$, and $K_1(\theta;\mathbb Z/q\mathbb Z)=0$ for every $q\ge2$.
\end{proposition}
\begin{proof}
Fix $u\in U(M_k(A))$. Apply Lemma~\ref{lem:localfd} to its matrix entries, with arbitrarily small error. We obtain a unital finite-dimensional subalgebra $F\subset B$ and $x\in M_k(F)$ with $\|x-\theta^{(k)}(u)\|$ arbitrarily small (the matrix norm is at most $k$ times the maximum entry norm). If this error is sufficiently small, both $x^*x$ and $xx^*$ are invertible in $M_k(F)$, and $v=x(x^*x)^{-1/2}\in U(M_k(F))$ is arbitrarily close to $\theta^{(k)}(u)$. Close unitaries represent the same $K_1$-class. Since $F$ is finite-dimensional, $K_1(F)=0$, and hence $[\theta^{(k)}(u)]=[v]=0$ in $K_1(B)$. This holds for every $k$ and every $u$.

For coefficients we use $K_*(A;\mathbb Z/q\mathbb Z)=K_*(A\otimes\mathcal O_{q+1})$, where $\otimes$ denotes the minimal tensor product (the only one, since $\mathcal O_{q+1}$ is nuclear), and $\theta\otimes\mathrm{id}$ is the induced map. Since $\mathcal O_{q+1}$ lies in the bootstrap class, the K\"unneth theorem gives a natural exact sequence $0\to K_1(A)\otimes\mathbb Z/q\mathbb Z\to K_1(A;\mathbb Z/q\mathbb Z)\to\operatorname{Tor}(K_0(A),\mathbb Z/q\mathbb Z)\to0$ \cite{Blackadar}. Let $u$ be a unitary in $M_k(A\otimes\mathcal O_{q+1})$. Approximate its entries by finite sums $\sum_ia_i\otimes o_i$ with $a_i\in A$ and $o_i\in\mathcal O_{q+1}$, and apply Lemma~\ref{lem:localfd} to the finite set of all $a_i$. This gives a unital finite-dimensional subalgebra $F\subset B$ such that $(\theta\otimes\mathrm{id})^{(k)}(u)$ is arbitrarily close to an element of $M_k(F\otimes\mathcal O_{q+1})$. As before, it is close to a unitary of $M_k(F\otimes\mathcal O_{q+1})$. Since $F\otimes\mathcal O_{q+1}$ is a finite direct sum of matrix algebras over $\mathcal O_{q+1}$ and $K_1(\mathcal O_{q+1})=0$, its class vanishes. This argument uses only Lemma~\ref{lem:localfd}, which does not require separability of $B$; so it applies to the nonseparable algebras $\Cu(\Gamma)$. For separable algebras the coefficient statement is also contained in the obstructions of \cite{CN2024}.
\end{proof}

\begin{theorem}\label{thm:zero}
If $\rdim(\Gamma)=0$, then $\Gamma$ is amenable and $K_1(\iota_\Gamma)=0$.
\end{theorem}
\begin{proof}
The reduced group C$^*$-algebra is separable and $\iota_\Gamma$ is unital, so Proposition~\ref{prop:qd} and Rosenberg's theorem imply amenability. Proposition~\ref{prop:kone} gives the second assertion.
\end{proof}

\subsection{The bilateral shift and the index pairing}
\begin{proposition}\label{prop:Zindex}
For $\Gamma=\mathbb Z$, the image in $K_1(\Cu(\mathbb Z))$ of the bilateral shift has infinite order. In particular, $\rdim(\mathbb Z)\ne0$.
\end{proposition}
\begin{proof}

Write $U\delta_n=\delta_{n+1}$ on $\ell^2(\mathbb Z)$, and let $P$ be the projection onto $\ell^2(\mathbb N_0)$. For each finite-propagation operator $T$, the two off-diagonal corners $PT(1-P)$ and $(1-P)TP$ have finite rank, since their matrix entries can be nonzero only within a finite distance of the cut between $-1$ and $0$. Thus $[P,T]$ is compact for every $T\in\Cu(\mathbb Z)$ by norm density. For $v\in U(M_k(\Cu(\mathbb Z)))$, the compression of $v$ by $P_k=1_k\otimes P$ to $P_k\ell^2(\mathbb Z)^k$ is Fredholm, with the compression of $v^*$ as parametrix modulo compacts. Homotopy invariance, stabilization and additivity of the Fredholm index define a homomorphism
\[
\operatorname{Ind}_P:K_1(\Cu(\mathbb Z))\longrightarrow\mathbb Z.
\]

The compression of $U$ is the unilateral shift, of index $-1$. Hence $\operatorname{Ind}_P([U])=-1$, so $[U]$ has infinite order and is in the image of $K_1(\iota_{\mathbb Z})$. Proposition~\ref{prop:kone} rules out dimension zero.
\end{proof}

\subsection{Characters and the Pimsner--Voiculescu boundary}
We now generalize the obstruction to amenable groups with arbitrary non-torsion abelianization. For the proof we allow coefficients larger than the usual group algebra, since the relevant finitely generated stages of $\Cu(\Gamma)$ have coefficient algebra $\ell^\infty(\Gamma)$, not $\ell^\infty(H)$.

\begin{lemma}\label{lem:pv}
Let a discrete group $H$ act on a unital C$^*$-algebra $D$ admitting an $H$-invariant tracial state. Suppose $q:H\twoheadrightarrow\mathbb Z$ is a homomorphism. There is a homomorphism
\[
I_q:K_1(D\rtimes_r H)\longrightarrow\mathbb R
\]
with $I_q([u_h])=q(h)$ for every $h\in H$, where $u_h$ is the canonical implementing unitary.
\end{lemma}
\begin{proof}
Let $N=\ker q$ and choose $t\in H$ with $q(t)=1$. Then $H=N\rtimes\langle t\rangle$. The decomposition $D\rtimes_r H\cong(D\rtimes_r N)\rtimes_r\mathbb Z$ used below holds for reduced crossed products by any such semidirect product, so no amenability of $H$ is needed. Set $B=D\rtimes_r N$; the invariant tracial state $\tau_D$ extends to a tracial state $\tau_B=\tau_D\circ E_N$ on $B$ by the canonical expectation $E_N:B\to D$. In fact, on finite Fourier sums the trace identity follows from the trace identity on $D$ and invariance under $N$. In particular $\tau_{B,*}([1_B])=1$. The action of $t$ induces an automorphism $\beta$ of $B$, and $D\rtimes_r H\cong B\rtimes_{\beta,r}\mathbb Z$. Let
\[
\partial:K_1(B\rtimes_\beta\mathbb Z)\longrightarrow K_0(B)
\]
be the Pimsner--Voiculescu boundary map \cite{PV1980,Blackadar}. With a fixed sign convention, $\partial([u_t])=[1_B]$ and $\partial$ vanishes on the image of $K_1(B)$. The first formula follows by naturality from the scalar inclusion $\mathbb C\to B$, or from the Toeplitz extension. Reverse the sign of $\partial$ if needed. For $h=nt^m$ with $n\in N$, we have $[u_h]=[u_n]+m[u_t]$ and therefore $\partial([u_h])=m[1_B]$. Taking $I_q=\tau_{B,*}\circ\partial$ proves the claim. The argument applies to nonseparable $B$ as well: the Toeplitz extension and its boundary map are defined without separability, or may be obtained by passage to invariant separable unital subalgebras and continuity of $K$-theory.
\end{proof}

\begin{theorem}\label{thm:characters}
Let $\Gamma$ be a countable amenable group. If $g\in\Gamma$ admits a homomorphism $\chi:\Gamma\to(\mathbb R,+)$ with $\chi(g)\ne0$, then the class of $\lambda_g$ has infinite order in $K_1(\Cu(\Gamma))$. In particular, non-torsion abelianization implies $\rdim(\Gamma)\ne0$.
\end{theorem}
\begin{proof}
Let $D=\ell^\infty(\Gamma)$ with its left translation action and choose a $\Gamma$-invariant mean on $D$, which exists by amenability. For every finitely generated subgroup $H\le\Gamma$ containing $g$, the image $\chi(H)$ is a finitely generated torsion-free abelian group, hence is isomorphic to $\mathbb Z^r$ for some $r\ge1$. There is a homomorphism $\chi(H)\to\mathbb Z$ nonzero on $\chi(g)$ (choose a suitable coordinate); after dividing by the positive generator of the image, the resulting map $q_H:H\twoheadrightarrow\mathbb Z$ is surjective and $q_H(g)\ne0$. Lemma~\ref{lem:pv}, applied to $D\rtimes_r H$ and the restricted invariant mean, shows that $[u_g]$ has infinite order in $K_1(D\rtimes_r H)$.

The crossed products $D\rtimes_r H$ form a directed family of unital subalgebras of $D\rtimes_r\Gamma$, via the canonical injective maps for subgroup inclusions, with dense union by the finite Fourier sum description. By continuity of $K$-theory \cite{Blackadar}, if any nonzero integer multiple of $[u_g]$ vanished in $K_1(D\rtimes_r\Gamma)$, it would already vanish at a finitely generated stage containing $g$. This contradicts the previous paragraph. Under the canonical identification with $\Cu(\Gamma)$, the class is the class of $\lambda_g$.

Finally, if $\Gamma_{\mathrm{ab}}$ is not torsion, choose a non-torsion element of it. Its image in $\Gamma_{\mathrm{ab}}\otimes_{\mathbb Z}\mathbb Q$ is nonzero; a rational linear functional nonzero on this image, followed by the inclusion $\mathbb Q\subset\mathbb R$, gives the required $\chi$ on $\Gamma$. Theorem~\ref{thm:zero} completes the proof.
\end{proof}

\begin{corollary}\label{cor:rational}
For amenable $\Gamma$, the homomorphism $g\mapsto[\lambda_g]$ induces an injective rational linear map
\[
\Gamma_{\mathrm{ab}}\otimes_{\mathbb Z}\mathbb Q
\longrightarrow K_1(\Cu(\Gamma))\otimes_{\mathbb Z}\mathbb Q.
\]
\end{corollary}
\begin{proof}
The equality $[uv]=[u]+[v]$ for unitaries in $K_1$ shows that the group homomorphism factors through the abelianization. Every rational tensor can be written as $\bar g\otimes(1/m)$ for an element $\bar g\in\Gamma_{\mathrm{ab}}$ and a positive integer $m$. Such a tensor is nonzero exactly when $\bar g$ is non-torsion. Theorem~\ref{thm:characters} shows that its image is nonzero as well.
\end{proof}

\section{Dimension one}

\begin{theorem}\label{thm:one}
If $\Gamma$ is exact, then $\rdim(\Gamma)\le1$. If $\Gamma$ is exact and non-amenable, then $\rdim(\Gamma)=1$.
\end{theorem}

\begin{proof}[First proof]
An exact group has Property A, so $\dimnuc(\Cu(\Gamma))\le1$ by \cite{GLZZ2026}, and $\rdim(\Gamma)\le1$ by Corollary~\ref{cor:upper}(a). If $\Gamma$ is also non-amenable, Theorem~\ref{thm:zero} gives $\rdim(\Gamma)\ge1$.
\end{proof}

\begin{proof}[Second proof of the upper bound]
This proof avoids \cite{GLZZ2026} and uses Proposition~\ref{prop:crossed}. If $\Gamma$ is finite, $\Cr(\Gamma)$ is finite-dimensional and $\rdim(\Gamma)=0$.

If $\Gamma$ is exact and non-amenable, R{\o}rdam and Sierakowski \cite{RS2012} construct a free, minimal, amenable action on the Cantor set $X$ such that $C(X)\rtimes_r\Gamma$ is a Kirchberg algebra satisfying the UCT. Kirchberg algebras have nuclear dimension one \cite{RSS, BGSW}.

If $\Gamma$ is amenable and infinite, Conley, Jackson, Kerr, Marks, Seward and Tucker-Drob \cite{CJKMST2018} show that a generic free minimal action of $\Gamma$ on the Cantor set $X$ has $\mathcal Z$-stable crossed product. This crossed product is separable, unital, nuclear, simple (the action is free and minimal \cite{AS1994}) and infinite-dimensional, so its nuclear dimension is at most one \cite{CETWW}.

In both cases Proposition~\ref{prop:crossed} gives $\rdim(\Gamma)\le1$.
\end{proof}

\begin{corollary}\label{cor:amenableone}
If $\Gamma$ is amenable with non-torsion abelianization, then $\rdim(\Gamma)=1$. In particular this holds for every countable amenable group with an infinite cyclic quotient, and for every non-torsion countable abelian group. Together with Theorem~\ref{thm:one}, every countable exact group with non-torsion abelianization has $\rdim(\Gamma)=1$.
\end{corollary}
\begin{proof}
Amenability gives exactness and hence the upper bound of Theorem~\ref{thm:one}. Theorem~\ref{thm:characters} gives the lower bound.
\end{proof}

\begin{corollary}\label{cor:abelian}
For a countable abelian group $A$,
\[
\rdim(A)=\begin{cases}0,&A\text{ is torsion},\\1,&A\text{ is not torsion}.\end{cases}
\]
\end{corollary}
\begin{proof}
A torsion abelian group is locally finite, so Corollary~\ref{cor:upper}(c) applies. Otherwise use Corollary~\ref{cor:amenableone}.
\end{proof}

\begin{example}\label{ex:one}
The invariant equals one for $\mathbb Z^n$ ($n\ge1$), $\mathbb Q$, $\bigoplus_{\mathbb N}\mathbb Z$, the discrete Heisenberg group, $\mathbb Z\wr\mathbb Z$, $\mathbb Z/2\mathbb Z\wr\mathbb Z$ and $BS(1,n)$ ($n\ge2$). The Heisenberg group has an infinite cyclic quotient through its abelianization; the wreath products have the quotient onto the acting copy of $\mathbb Z$, and the Baumslag--Solitar group has the exponent-sum homomorphism onto $\mathbb Z$. More generally, every infinite finitely generated nilpotent group has infinite abelianization: if its abelianization were finite, the successive lower central quotients, generated by iterated commutators and quotients of tensor powers of the abelianization, would be finite, so the group would be finite. Thus all infinite finitely generated nilpotent groups have relative nuclear dimension one.
\end{example}

\section{Dihedral groups and permanence of dimension zero}\label{sec:dihedral}

The abelian calculation might suggest that dimension zero detects local finiteness or at least the absence of an infinite cyclic subgroup. The next theorem rules out both possibilities. Its proof uses a classical AF crossed-product result rather than a one-colour approximation constructed directly on the group.

\begin{theorem}\label{thm:dihedral}
Let $D_\infty=\langle t,s\mid s^2=e,\ sts=t^{-1}\rangle$. Then $\rdim(D_\infty)=0$, while $\dimnuc(\Cu(D_\infty))=1$. In fact, $\Cr(D_\infty)$ is contained unitally in a separable AF subalgebra of $\Cu(D_\infty)$.
\end{theorem}
\begin{proof}
Let $X=\varprojlim_k\mathbb Z/2^k\mathbb Z$ be the additive group of $2$-adic integers, viewed as a Cantor space. Define homeomorphisms $T(x)=x+1$ and $S(x)=-x$. They satisfy $S^2=\mathrm{id}$ and $STS=T^{-1}$, defining an action of $D_\infty$ by $t\mapsto T$ and $s\mapsto S$. The homeomorphism $T$ is minimal: its orbit meets every basic clopen cylinder because the integers map onto each finite quotient. The involution $S$ fixes $0\in X$.

Bratteli, Evans and Kishimoto \cite[Theorem~3.5]{BEK1993} prove that for a minimal homeomorphism of a totally disconnected compact metrizable space and a reversing involution, the crossed product by $D_\infty$ is AF when the involution or its product with the minimal homeomorphism has a fixed point. These hypotheses hold here. Since $D_\infty$ is amenable, the full and reduced crossed products agree. Thus
\[
B=C(X)\rtimes_r D_\infty
\]
is a separable unital AF algebra. The action need not be free, and no freeness is required in Proposition~\ref{prop:crossed}. Apply that proposition with $x_0=0$. Since the $D_\infty$-orbit of $0$ is dense (already its $\langle t\rangle$-orbit is dense), the coefficient map $C(X)\to\ell^\infty(D_\infty)$ is injective. Functoriality of reduced crossed products for injective equivariant maps then makes $B\to\Cu(D_\infty)$ injective. Its image is a unital AF subalgebra containing $\Cr(D_\infty)$, and Lemma~\ref{lem:factor} gives $\rdim(D_\infty)=0$. The group is exact (it is amenable) but not locally finite (it contains $\langle t\rangle\cong\mathbb Z$). By (F3) its uniform Roe algebra has nuclear dimension one.
\end{proof}

\begin{corollary}\label{cor:nonmonotone}
The invariant $\rdim$ is not monotone under passage to subgroups, not invariant under finite-index inclusions, and not invariant under coarse equivalence of groups endowed with proper invariant metrics.
\end{corollary}
\begin{proof}
The inclusion $\langle t\rangle\cong\mathbb Z\le D_\infty$ has index two. Corollary~\ref{cor:abelian} and Theorem~\ref{thm:dihedral} give $1=\rdim(\mathbb Z)>\rdim(D_\infty)=0$. The two groups are finitely generated and the inclusion of a finite-index subgroup is a quasi-isometry for word metrics (equivalently for right-invariant word metrics), hence a coarse equivalence. The corresponding proper invariant metrics induce the same coarse structures.
\end{proof}

\begin{remark}
There is no contradiction with Proposition~\ref{prop:Zindex}. The block representation of the proof of Proposition~\ref{prop:permanence}(ii), applied to $H=\langle t\rangle$ with coset representatives $\{e,s\}$, is a unital $*$-homomorphism $\Cu(\langle t\rangle)\to\Cu(D_\infty)$ sending the bilateral shift $U$ to $\lambda_t$. So the infinite-order class $[U]\in K_1(\Cu(\mathbb Z))$ is sent to $[\lambda_t]$, which vanishes in $K_1(\Cu(D_\infty))$ because the regular inclusion of $D_\infty$ factors through an AF algebra. The geometric reason is that the second coset runs backwards. For the right-invariant word metric one has $|(t^ns)(t^m)^{-1}|=|t^{n+m}s|$, so $t^ns$ lies at bounded distance from $t^{-n}$. Hence $\lambda_t$ moves the two cosets in opposite directions, and for the projection $P$ onto $\{t^n,\,t^{-n}s:n\ge0\}$, which has compact commutator with every element of $\Cu(D_\infty)$, the compression of $\lambda_t$ has Fredholm index $-1+1=0$. Also $D_{\infty,\mathrm{ab}}\cong\mathbb Z/2\mathbb Z\oplus\mathbb Z/2\mathbb Z$, consistent with Theorem~\ref{thm:characters}.
\end{remark}

\begin{proposition}\label{prop:permanence}
\begin{enumerate}[label=(\roman*)]
\item If $\rdim(\Gamma)=\rdim(\Lambda)=0$, then $\rdim(\Gamma\times\Lambda)=0$.
\item If $\Gamma=\bigcup_{j\ge1}\Gamma_j$ is an increasing union of subgroups satisfying $\rdim(\Gamma_j)=0$, then $\rdim(\Gamma)=0$.
\end{enumerate}
\end{proposition}
\begin{proof}
(i) Under $\ell^2(\Gamma\times\Lambda)=\ell^2(\Gamma)\otimes\ell^2(\Lambda)$, the spatial tensor product $\Cu(\Gamma)\otimes_{\min}\Cu(\Lambda)$ embeds unitally into $\Cu(\Gamma\times\Lambda)$: tensor products of multiplication operators are multiplication operators on the product group, and the tensor products of the regular unitaries are regular unitaries of the product group. The algebra $\Cr(\Gamma\times\Lambda)$ identifies with $\Cr(\Gamma)\otimes_{\min}\Cr(\Lambda)$. Given a finite set and an error tolerance, approximate its elements by finite sums of elementary tensors. Lemma~\ref{lem:localfd} provides unital finite-dimensional subalgebras $F_\Gamma\subset\Cu(\Gamma)$ and $F_\Lambda\subset\Cu(\Lambda)$ approximating all factors of these sums to any prescribed accuracy. The elementary estimate
\[
\|a\otimes b-x\otimes y\|\le\|a-x\|\|b\|+\|x\|\|b-y\|
\]
shows that the tensor product $F_\Gamma\otimes F_\Lambda$ approximates the initial finite set, with arbitrarily small error. It is a unital finite-dimensional subalgebra of the product's uniform Roe algebra. Apply Lemma~\ref{lem:localfd} again.

(ii) We first describe an embedding for any subgroup $H\le\Gamma$. Choose representatives $R$ for the right cosets $Hr$ of $H$ in $\Gamma$. Then $\ell^2(\Gamma)\cong\bigoplus_{r\in R}\ell^2(H)$ via $\delta_{hr}\mapsto\delta_h$ in the $r$-summand. The repeated block representation $T\mapsto\bigoplus_{r\in R}T$ of $\Cu(H)$ maps a diagonal function $f\in\ell^\infty(H)$ to multiplication by the bounded function $hr\mapsto f(h)$ on $\Gamma$, and it maps $\lambda_h^H$ to $\lambda_h^\Gamma$ for every $h\in H$. Hence it is a unital injective $*$-homomorphism $\Cu(H)\to\Cu(\Gamma)$ extending the canonical group algebra inclusion. Now every finite set in $\Cr(\Gamma)$ can be approximated arbitrarily well by elements of $\Cr(\Gamma_j)$ for one common $j$, since finite group-algebra sums are dense. Use the local finite-dimensional approximation for $\Gamma_j$ and embed the resulting algebra into $\Cu(\Gamma)$ by the preceding construction. Lemma~\ref{lem:localfd} gives the result.
\end{proof}

\begin{corollary}\label{cor:products}
Let $L$ be any countable locally finite group and $n\ge1$. Then $\rdim(D_\infty^n\times L)=0$. In particular, for every $n\ge1$ there are quasi-isometric finitely generated amenable groups $D_\infty^n$ and $\mathbb Z^n$ of asymptotic dimension $n$ with distinct relative nuclear dimensions.
\end{corollary}
\begin{proof}
Use Corollary~\ref{cor:upper}(c), Theorem~\ref{thm:dihedral} and Proposition~\ref{prop:permanence}. The subgroup $\mathbb Z^n\le D_\infty^n$ has finite index $2^n$, and $\asdim\mathbb Z^n=n$; hence the groups have the same asymptotic dimension. Corollary~\ref{cor:abelian} gives the contrasting value one.
\end{proof}

\subsection{Holonomy obstructions for virtually abelian groups}
The $K_1$-obstruction also sees classes of higher degree. Recall that $K_*(\Cr(\mathbb Z^n))\cong K^*(\mathbb T^n)$ is isomorphic, as a $\mathbb Z/2\mathbb Z$-graded ring, to the exterior algebra $\Lambda^*(\mathbb Z^n)$, with the classes $[\lambda_{e_i}]$ of the coordinate unitaries corresponding to the standard basis vectors $e_i$ of $\Lambda^1(\mathbb Z^n)$; the even part is $K_0$ and the odd part is $K_1$. For $I=\{i_1<\dots<i_r\}\subset\{1,\dots,n\}$ write $e_I=e_{i_1}\wedge\cdots\wedge e_{i_r}$, with $e_\emptyset=[1]$. The automorphism $\alpha_f$ of $\Cr(\mathbb Z^n)$ induced by $f\in\mathrm{GL}(n,\mathbb Z)$ acts on $\Lambda^*(\mathbb Z^n)$ by $\Lambda^*(f)$, and on the top degree by multiplication with $\det f$. For a subset $I$ let $H_I=\bigoplus_{i\in I}\mathbb Ze_i$.

\begin{theorem}\label{thm:injective}
For every $n\ge1$, the inclusion $\iota_{\mathbb Z^n}\colon\Cr(\mathbb Z^n)\hookrightarrow\Cu(\mathbb Z^n)$ induces injective maps $K_j(\Cr(\mathbb Z^n))\to K_j(\Cu(\mathbb Z^n))$ for $j=0,1$.
\end{theorem}

\begin{proof}
Let $H=\mathbb Z^n$ and $D=\ell^\infty(H)$, and fix an invariant mean $\mu$ on $D$. Let $x=\sum_Ia_Ie_I\ne0$ be a class of fixed parity and choose $I_0$ of maximal cardinality with $a_{I_0}\ne0$. Write $H=H_{I_0^c}\oplus H_{I_0}$, so that
\[
\Cu(H)=D\rtimes_rH=\big(D\rtimes_rH_{I_0^c}\big)\rtimes H_{I_0},\qquad \Cr(H)=\Cr(H_{I_0^c})\rtimes H_{I_0},
\]
the second with trivial action, and the inclusion $\Cr(H_{I_0^c})\to B=D\rtimes_rH_{I_0^c}$ induced by $\mathbb C1\subset D$ is compatible with these decompositions. Crossing by the coordinates of $I_0$ one at a time and composing the Pimsner--Voiculescu boundary maps \cite{PV1980} gives homomorphisms $\partial_{I_0}$ on both sides, related by naturality. For the trivial actions on the group algebra side, each boundary map is, up to sign, contraction with the corresponding coordinate. So $\partial_{I_0}(e_I)=0$ unless $I_0\subseteq I$, and $\partial_{I_0}(e_{I_0})=\pm[1]$. By maximality of $|I_0|$, the only term of $x$ with $I\supseteq I_0$ is $a_{I_0}e_{I_0}$, so $\partial_{I_0}(x)=\pm a_{I_0}[1]$ and
\[
\partial_{I_0}\big(K_*(\iota_{\mathbb Z^n})(x)\big)=\pm a_{I_0}[1_B]\in K_0(B).
\]

The state $\tau=\mu\circ E$, with $E\colon B\to D$ the canonical conditional expectation, is tracial because $\mu$ is invariant, and $\tau_*([1_B])=1$. So $\tau_*\partial_{I_0}(K_*(\iota_{\mathbb Z^n})(x))=\pm a_{I_0}\ne0$. If $I_0=\emptyset$, then $x=a_\emptyset[1]$ and the same trace on $\Cu(H)$ gives the conclusion directly. As in Lemma~\ref{lem:pv}, separability plays no role.
\end{proof}

\begin{remark}
The $K$-groups of $\Cu(\mathbb Z^n)$ themselves have been computed. Kato, Kishimoto and Tsutaya treat $n\le2$ \cite{KKT2023}, and Fanning and Aktas compute $K_*(\Cu(G,\mathcal E))$ for countable abelian groups $G$ with the coarse structure generated by a finite-rank free abelian subgroup, by a spectral sequence in group homology \cite{FA2026}. For instance, with $\mathfrak c=2^{\aleph_0}$, they obtain $K_0(\Cu(\mathbb Z))\cong\mathbb Q^{(\mathfrak c)}$ and $K_1(\Cu(\mathbb Z))\cong\mathbb Z$, the latter consistent with Proposition~\ref{prop:Zindex}, and $K_0(\Cu(\mathbb Z^2))\cong\mathbb Q^{(\mathfrak c)}\oplus\mathbb Z$, $K_1(\Cu(\mathbb Z^2))\cong\mathbb Q^{(\mathfrak c)}$. These computations are not stated in terms of the map induced by the regular inclusion, which is what Theorem~\ref{thm:injective} describes: it shows that the finitely generated group $K_*(\Cr(\mathbb Z^n))\cong\Lambda^*(\mathbb Z^n)$ embeds, and the iterated boundary maps composed with an invariant trace give explicit coordinates on its image.
\end{remark}

\begin{lemma}\label{lem:block}
Let $\Gamma$ be a countable group with a normal subgroup $H$ of finite index $m$, and let $r_1,\dots,r_m$ be representatives of the left cosets of $H$. The unitary $W\colon\ell^2(\Gamma)\to\ell^2(H)^m$, $\delta_{r_jh}\mapsto\delta_h$ in the $j$-th summand, induces a $*$-isomorphism $\rho=\mathrm{Ad}\,W\colon\Cu(\Gamma)\to M_m(\Cu(H))$. It maps $\Cr(\Gamma)$ into $M_m(\Cr(H))$, so $\rho\circ\iota_\Gamma=\iota_H^{(m)}\circ\rho_0$ for a $*$-homomorphism $\rho_0\colon\Cr(\Gamma)\to M_m(\Cr(H))$. If $\kappa\colon\Cr(H)\to\Cr(\Gamma)$ is the canonical inclusion and $f_j(h)=r_j^{-1}hr_j$, then $\rho_0\circ\kappa=\mathrm{diag}(\alpha_{f_1},\dots,\alpha_{f_m})$.
\end{lemma}

\begin{proof}
Fix a proper right-invariant metric on $\Gamma$ with length function $|\cdot|$, and give $H$ the restricted metric; this is a proper right-invariant metric on $H$, so it defines $\Cu(H)$. Put $M=\max_j|r_j|$. If $T\in\Cu(\Gamma)$ has propagation at most $R$, the $(j,k)$ block of $WTW^*$ has nonzero entries only at pairs $(h,h')$ with $|r_jhh'^{-1}r_k^{-1}|\le R$, and then $|hh'^{-1}|\le R+2M$. Conversely, if all blocks of a matrix $(T_{jk})$ have propagation at most $R$ in $H$, then $W^*(T_{jk})W$ has entries only at pairs $(r_jh,r_kh')$ with $|hh'^{-1}|\le R$, and $d(r_jh,r_kh')=|r_jhh'^{-1}r_k^{-1}|\le R+2M$. Multiplication operators correspond to diagonal matrices of multiplication operators. So $\mathrm{Ad}\,W$ maps the finite-propagation operators of $\Gamma$ onto matrices of finite-propagation operators of $H$, and $\rho$ is a $*$-isomorphism. For $g\in\Gamma$ write $gr_j=r_{p(j)}k_j$ with $k_j\in H$. Then $\lambda_g\delta_{r_jh}=\delta_{r_{p(j)}k_jh}$, so $\rho(\lambda_g)$ is the permutation matrix of $p$ with entries $\lambda^H_{k_j}$, which lies in $M_m(\Cr(H))$. For $g=h\in H$ one has $hr_j=r_jf_j(h)$, which gives the last formula.
\end{proof}

\begin{theorem}\label{thm:holonomy}
Let $\Gamma$ be a countable group with a normal subgroup $H\cong\mathbb Z^n$ of finite index, and let $F\subset\mathrm{GL}(n,\mathbb Z)$ be the holonomy group, the image of the conjugation action of $\Gamma$ on $H$. Then the image of $K_1(\iota_\Gamma)$ is a free abelian group of rank
\[
\sum_{k\ge0}\dim_{\mathbb Q}\Big(\bigwedge\nolimits^{2k+1}\mathbb Q^n\Big)^F.
\]

Moreover, let $L_1=(\Lambda^{\mathrm{odd}}\mathbb Z^n)^F$ and $m=[\Gamma:H]$. Under the identifications of Lemma~\ref{lem:block}, the image of $K_1(\iota_\Gamma)$ corresponds to a subgroup $M_1$ with $mL_1\subseteq M_1\subseteq L_1$, so its index in $L_1$ divides $m^{\operatorname{rank}L_1}$. In particular $K_1(\iota_\Gamma)=0$ if and only if $(\Lambda^{2k+1}\mathbb Q^n)^F=0$ for all $k\ge0$, and if $(\Lambda^{2k+1}\mathbb Q^n)^F\ne0$ for some $k$, then $\rdim(\Gamma)=1$.
\end{theorem}

\begin{proof}
WE use the notation of Lemma~\ref{lem:block} in this proof. We first show that the image of $K_*(\rho_0)$ consists of $F$-invariant classes, under $K_*(M_m(\Cr(H)))\cong K_*(\Cr(H))$. Let $g\in\Gamma$ and let $R_g\delta_x=\delta_{xg^{-1}}$ be the right translation, which commutes with $\lambda(\Gamma)$. Put $\varphi_g(h)=ghg^{-1}$ and write $r_jg^{-1}=r_{q(j)}k'_j$ with $k'_j\in H$. Then $R_g\delta_{r_jh}=\delta_{r_{q(j)}k'_j\varphi_g(h)}$, so
\[
WR_gW^*=V_g\,U_{\varphi_g}^{(m)},
\]
where $U_{\varphi_g}\delta_h=\delta_{\varphi_g(h)}$ implements $\alpha_{\varphi_g}$ on $\Cr(H)$, $U_{\varphi_g}^{(m)}=\mathrm{diag}(U_{\varphi_g},\dots,U_{\varphi_g})$ acts on $\ell^2(H)^m$, and $V_g$ is the unitary permutation matrix of $q$ with entries $\lambda^H_{k'_j}$, an element of $M_m(\Cr(H))$. Since $R_g$ commutes with $\lambda(\Gamma)$,
\[
\rho_0(a)=V_g\,\alpha_{\varphi_g}^{(m)}\big(\rho_0(a)\big)\,V_g^*\qquad(a\in\Cr(\Gamma)).
\]

Inner automorphisms act trivially on $K$-theory, so $K_*(\rho_0)=K_*(\alpha_{\varphi_g})\circ K_*(\rho_0)$. As $g$ runs through $\Gamma$, $\varphi_g$ runs through $F$.

Conversely, if $c\in K_1(\Cr(H))$ is $F$-invariant, Lemma~\ref{lem:block} gives
\[
K_1(\rho_0)K_1(\kappa)(c)=\sum_{j=1}^mK_1(\alpha_{f_j})(c)=mc.
\]

So the image $M_1$ of $K_1(\rho_0)$ satisfies $mL_1\subseteq M_1\subseteq L_1$ inside the free abelian group $\Lambda^{\mathrm{odd}}(\mathbb Z^n)$. Since $[L_1:mL_1]=m^{\operatorname{rank}L_1}$, the index of $M_1$ in $L_1$ divides $m^{\operatorname{rank}L_1}$, and the rank of $M_1$ is $\dim_{\mathbb Q}(\Lambda^{\mathrm{odd}}\mathbb Q^n)^F$, the sum in the statement.

Finally, $K_1(\rho)$ is an isomorphism and $K_1(\iota_H^{(m)})$ is injective by Theorem~\ref{thm:injective} and stability of $K$-theory. From $\rho\circ\iota_\Gamma=\iota_H^{(m)}\circ\rho_0$, the image of $K_1(\iota_\Gamma)$ is isomorphic to the image of $K_1(\rho_0)$. If it is nonzero, $\rdim(\Gamma)\ne0$ by Proposition~\ref{prop:kone}; $\Gamma$ is virtually abelian, hence amenable, and $\rdim(\Gamma)\le1$ by Theorem~\ref{thm:one}.
\end{proof}

The invariance of the image of $K_*(\rho_0)$, the transfer formula and Theorem~\ref{thm:injective} hold in both degrees, so the same argument shows that the image of $K_0(\iota_\Gamma)$ is free abelian of rank $\sum_k\dim_{\mathbb Q}(\Lambda^{2k}\mathbb Q^n)^F$, and corresponds to a subgroup $M_0$ with $mL_0\subseteq M_0\subseteq L_0$, where $L_0=(\Lambda^{\mathrm{even}}\mathbb Z^n)^F$. The first exterior degree recovers the character obstruction of Theorem~\ref{thm:characters} for virtually abelian groups: since $F$ is finite, $\Gamma_{\mathrm{ab}}\otimes\mathbb Q\cong(\mathbb Q^n)_F\cong(\mathbb Q^n)^F$.

\begin{corollary}\label{cor:determinant}
Under the hypotheses of Theorem~\ref{thm:holonomy},
\[
\begin{aligned}
\operatorname{rank}\operatorname{im}K_1(\iota_\Gamma)&=\frac1{2|F|}\sum_{f\in F}\big(\det(I+f)-\det(I-f)\big),\\
\operatorname{rank}\operatorname{im}K_0(\iota_\Gamma)&=\frac1{2|F|}\sum_{f\in F}\big(\det(I+f)+\det(I-f)\big).
\end{aligned}
\]
\end{corollary}

\begin{proof}
For a finite group, the averaging operator $\frac1{|F|}\sum_{f\in F}\Lambda^j(f)$ is a projection of $\Lambda^j\mathbb Q^n$ onto $(\Lambda^j\mathbb Q^n)^F$, so $\dim_{\mathbb Q}(\Lambda^j\mathbb Q^n)^F=\frac1{|F|}\sum_{f\in F}\operatorname{Tr}\Lambda^j(f)$. The identity $\det(I+tf)=\sum_{j=0}^nt^j\operatorname{Tr}\Lambda^j(f)$ at $t=1$ and $t=-1$ gives
\[
\sum_{j\ \mathrm{odd}}\operatorname{Tr}\Lambda^j(f)=\tfrac12\big(\det(I+f)-\det(I-f)\big),\qquad
\sum_{j\ \mathrm{even}}\operatorname{Tr}\Lambda^j(f)=\tfrac12\big(\det(I+f)+\det(I-f)\big).
\]

Average over $F$ and apply Theorem~\ref{thm:holonomy} and the remark following it.
\end{proof}

For $D_\infty^n$, with $F$ all diagonal sign matrices, $\det(I+f)$ is $2^n$ for $f=I$ and $0$ otherwise, and $\det(I-f)$ is $2^n$ for $f=-I$ and $0$ otherwise. So the image of $K_1(\iota_{D_\infty^n})$ is zero and the image of $K_0(\iota_{D_\infty^n})$ has rank one.

The integral image can be identified. Let $S\cong(\mathbb Z/2\mathbb Z)^n$ be the subgroup of $D_\infty^n$ generated by the $n$ reflections and $p_S=2^{-n}\sum_{s\in S}\lambda_s$, a projection in $\Cr(D_\infty^n)$. With the elements of $S$ as coset representatives for $H=\mathbb Z^n$, Lemma~\ref{lem:block} sends each $\lambda_s$ to a permutation matrix with scalar entries, and $\rho_0(p_S)$ is the constant rank-one projection $2^{-n}(1)_{s,t\in S}$ in $M_{2^n}(\mathbb C)\subset M_{2^n}(\Cr(H))$. So $K_0(\rho_0)[p_S]=[1]$ and $K_0(\rho_0)[1]=2^n[1]$. Here $L_0=\mathbb Z[1]$, since only degree zero carries invariants, so $M_0=L_0$ and the image of $K_0(\iota_{D_\infty^n})$ is infinite cyclic, generated by $[\iota(p_S)]$, with $[1_{\Cu(D_\infty^n)}]=2^n[\iota(p_S)]$. The unit therefore spans a subgroup of index $2^n$. A trace confirms that the unit is not a generator. Since $D_\infty^n$ is amenable, an invariant mean $\mu$ on $\ell^\infty(D_\infty^n)$ gives a tracial state $\tilde\tau=\mu\circ E$ on $\Cu(D_\infty^n)$, where $E$ is the canonical conditional expectation onto $\ell^\infty(D_\infty^n)$. It restricts to the canonical trace on $\Cr(D_\infty^n)$, so $\tilde\tau_*[\iota(p_S)]=2^{-n}$ and $\tilde\tau_*[1]=1$. So $[\iota(p_S)]$ is not an integer multiple of $[1]$.

\begin{corollary}\label{cor:bieberbach}
Let $\Gamma$ be a nontrivial finitely generated torsion-free virtually abelian group. Then the images of $K_0(\iota_\Gamma)$ and $K_1(\iota_\Gamma)$ have the same rank, which is at least one, and $\rdim(\Gamma)=1$.
\end{corollary}

\begin{proof}
The normal core of a finite-index abelian subgroup is a finite-index normal subgroup $H$, which is finitely generated, torsion-free and abelian, so $H\cong\mathbb Z^n$; and $n\ge1$ because $\Gamma$ is infinite. Let $g\in\Gamma$ and let $f$ be the automorphism $h\mapsto g^{-1}hg$ of $H$. We claim that $f$ has eigenvalue $1$. Choose $k\ge1$ with $g^k\in H$. Then $f(g^k)=g^k$. If $f$ had no eigenvalue $1$, then $f-I$ would be injective on $H$, so $g^k=e$, and $g=e$ since $\Gamma$ is torsion-free. Then $f=I$, which has eigenvalue $1$, a contradiction. So $\det(I-f)=0$ for every $f$ in the holonomy group $F$, and Corollary~\ref{cor:determinant} gives $\operatorname{rank}\operatorname{im}K_0(\iota_\Gamma)=\operatorname{rank}\operatorname{im}K_1(\iota_\Gamma)$. The rank in $K_0$ is at least one, since $\Lambda^0\mathbb Q^n$ is always invariant. So $K_1(\iota_\Gamma)\ne0$, and Theorem~\ref{thm:holonomy} gives $\rdim(\Gamma)=1$.
\end{proof}

In particular every Bieberbach group of positive rank has relative nuclear dimension one, and among nontrivial finitely generated virtually abelian groups, relative nuclear dimension zero requires nontrivial torsion, as in $D_\infty^n$.

\begin{corollary}\label{cor:orientation}
Let $\Gamma$ have a normal subgroup $H\cong\mathbb Z^n$ of finite index, with $n$ odd, such that conjugation by every element of $\Gamma$ acts on $H$ with determinant one. Then $\rdim(\Gamma)=1$.
\end{corollary}

\begin{proof}
The class $e_1\wedge\cdots\wedge e_n\in\Lambda^n(\mathbb Z^n)$ lies in odd degree and is fixed by $F\subset\mathrm{SL}(n,\mathbb Z)$. Apply Theorem~\ref{thm:holonomy}.
\end{proof}

The determinant condition cannot be dropped. For $D_\infty$ ($n=1$, holonomy $\{\pm1\}$) and for $D_\infty^n$ (holonomy all diagonal sign matrices), a diagonal sign matrix multiplies $e_I$ by the product of its signs over $I$, so no nonzero exterior class of positive degree is invariant; this is consistent with Theorem~\ref{thm:dihedral} and Corollary~\ref{cor:products}.

\begin{example}\label{ex:orientation}
Let $V=\{\mathrm{diag}(\epsilon_1,\epsilon_2,\epsilon_3):\epsilon_i=\pm1,\ \epsilon_1\epsilon_2\epsilon_3=1\}\cong(\mathbb Z/2\mathbb Z)^2$, generated by $A=\mathrm{diag}(1,-1,-1)$ and $B=\mathrm{diag}(-1,1,-1)$, and let $\Gamma_0=\mathbb Z^3\rtimes V$. Writing $D_\infty^3=\mathbb Z^3\rtimes\{\pm1\}^3$ with the diagonal sign matrices, $\Gamma_0$ is the subgroup of index two given by the condition $\det=1$. The abelianization of a split extension is $(\Gamma_0)_{\mathrm{ab}}\cong(\mathbb Z^3)_V\oplus V$, and $(1-A)\mathbb Z^3+(1-B)\mathbb Z^3=2\mathbb Z\oplus2\mathbb Z\oplus2\mathbb Z$, so $(\Gamma_0)_{\mathrm{ab}}\cong(\mathbb Z/2\mathbb Z)^5$. Corollary~\ref{cor:orientation} gives $\rdim(\Gamma_0)=1$, while $\rdim(D_\infty^3)=0$ by Corollary~\ref{cor:products}. In Corollary~\ref{cor:determinant}, $\det(I+f)=8$ for $f=I$ and $0$ for $f=A,B,AB$, while $\det(I-f)=0$ for all $f\in V$; so both images have rank one. Under the identifications of Lemma~\ref{lem:block}, the $K_1$-image lies between $4\,\mathbb Z\,e_1\wedge e_2\wedge e_3$ and $\mathbb Z\,e_1\wedge e_2\wedge e_3$; rationally it is the line spanned by $e_1\wedge e_2\wedge e_3$. We have not determined its index.

The group $\Gamma_4=\Gamma_0\times D_\infty$ has the normal subgroup $\mathbb Z^4$ of index $8$, with holonomy $V\times\{\pm1\}$, and $(\Gamma_4)_{\mathrm{ab}}\cong(\mathbb Z/2\mathbb Z)^7$. The odd class $e_1\wedge e_2\wedge e_3$ is invariant, so $\rdim(\Gamma_4)=1$ by Theorem~\ref{thm:holonomy}, while $\rdim(D_\infty^4)=0$. Here the rank of the free abelian subgroup is even, so Corollary~\ref{cor:orientation} does not apply.
\end{example}

\begin{corollary}\label{cor:torsionab}
Amenability together with torsion abelianization does not imply $\rdim(\Gamma)=0$. Passing to a subgroup of index two can raise the relative nuclear dimension.
\end{corollary}

\begin{proof}
The group $\Gamma_0$ of Example~\ref{ex:orientation} is amenable, has abelianization $(\mathbb Z/2\mathbb Z)^5$ and satisfies $\rdim(\Gamma_0)=1$. It has index two in $D_\infty^3$, for which $\rdim(D_\infty^3)=0$.
\end{proof}

\begin{remark}
For virtually abelian groups the $K_1$-obstruction passes to finite-index subgroups. If $\Lambda\le\Gamma$ has finite index and $H\cong\mathbb Z^n$ is normal of finite index in $\Gamma$, then $H\cap\Lambda\cong\mathbb Z^n$ is normal of finite index in $\Lambda$. Its holonomy group, acting on $(H\cap\Lambda)\otimes\mathbb Q=\mathbb Q^n$, is a subgroup of $F$, so it has at least as many invariants in each exterior power. So a finite extension of a group with nonzero $K_1$-obstruction may lose the obstruction, as for $\mathbb Z\le D_\infty$ and $\Gamma_0\le D_\infty^3$, but a finite-index subgroup of a group with nonzero obstruction keeps it.
\end{remark}

\begin{proof}[Proof of the Main Theorem]
Part (i) follows from Proposition~\ref{prop:exact} and Theorem~\ref{thm:one}. Part (ii) follows from Theorems~\ref{thm:zero} and \ref{thm:characters}. Part (iii) is Corollary~\ref{cor:upper}(c), Theorem~\ref{thm:one}, Corollary~\ref{cor:amenableone} and Theorem~\ref{thm:holonomy}. Part (iv) is Corollary~\ref{cor:abelian}; part (v) is Theorem~\ref{thm:dihedral} and Corollary~\ref{cor:nonmonotone}. Part (vi) follows from Proposition~\ref{prop:permanence} and Corollary~\ref{cor:products}, part (vii) is Example~\ref{ex:orientation}, and part (viii) is Theorems~\ref{thm:injective} and \ref{thm:holonomy} and Corollary~\ref{cor:bieberbach}.
\end{proof}

\section{Comparison and questions}\label{sec:questions}

\begin{center}
\small
\renewcommand{\arraystretch}{1.25}
\begin{tabular}{lccc}
\hline
$\Gamma$ & $\dimnuc(\Cr(\Gamma))$ & $\dimnuc(\Cu(\Gamma))$ & $\rdim(\Gamma)$\\
\hline
locally finite & $0$ & $0$ & $0$\\
$\mathbb Z^n$, $n\ge1$ & $n$ & $1$ & $1$\\
$\mathbb Z\wr\mathbb Z$ & $\infty$ & $1$ & $1$\\
$\bigoplus_{\mathbb N}\mathbb Z$ & $\infty$ & $1$ & $1$\\
$D_\infty$ & $1$ & $1$ & $0$\\
$D_\infty^n$, $n\ge2$ & $n$ & $1$ & $0$\\
$\Gamma_0=\mathbb Z^3\rtimes V$ (Example~\ref{ex:orientation}) & $3$ & $1$ & $1$\\
$\Gamma_4=\Gamma_0\times D_\infty$ & $4$ & $1$ & $1$\\
$\mathbb F_n$, $n\ge2$ & $\infty$ & $1$ & $1$\\
$\mathrm{SL}(3,\mathbb Z)$ & $\infty$ & $1$ & $1$\\
non-exact groups & $\infty$ & $\infty$ & $\infty$\\
\hline
\end{tabular}
\end{center}

The third column follows from \cite{GLZZ2026,LW2018}. The value for $\mathbb Z\wr\mathbb Z$ in the second column follows from Eckhardt and Wu \cite[Theorem~C]{EW2026}. The dual of $\bigoplus_{\mathbb N}\mathbb Z$ is the infinite product of circles, which has infinite covering dimension. The infinite dihedral group algebra is the interval algebra of continuous $M_2$-valued functions with diagonal values at the two endpoints, by Fourier transform and the reflection action on the circle. It has nuclear dimension one: the decomposition rank formula for subhomogeneous algebras \cite{Winter2004}, together with nuclear dimension being bounded above by decomposition rank, gives an upper bound of one and its ideal $C_0((0,1),M_2)$ rules out dimension zero (which would make it AF). The remaining entries of the second column for virtually abelian groups follow from one argument. Each of $D_\infty^n$, $\Gamma_0$ and $\Gamma_4$ is of the form $\mathbb Z^n\rtimes V$ with $V$ a finite group of diagonal sign matrices, and Fourier transform gives $\Cr(\mathbb Z^n\rtimes V)\cong C(\mathbb T^n)\rtimes V$, a subhomogeneous algebra. By the Mackey machine, every irreducible representation is induced from a point $z\in\mathbb T^n$ together with an irreducible representation of its stabilizer $V_z$, and its dimension is $[V:V_z]\dim\sigma$. So the primitive ideal space is a finite union of pieces, one for each subgroup $S\le V$ and each irreducible representation of $S$, and the piece for $S$ is homeomorphic to the quotient of $\{z\in\mathbb T^n:V_z=S\}$ by the finite group $V$. These sets are locally closed in $\mathbb T^n$, and a quotient of a subset of $\mathbb T^n$ by a finite group has covering dimension at most $n$. Each stratum of irreducible representations of a fixed dimension is a finite union of such pieces, which are locally closed and hence $F_\sigma$ in the metrizable stratum, so it has covering dimension at most $n$ by the countable sum theorem. The decomposition rank formula for subhomogeneous algebras \cite{Winter2004} gives $\dimnuc\le n$. For the lower bound, $V$ acts freely on the dense open set $W$ of points with trivial stabilizer, so $C_0(W)\rtimes V$ is an ideal which is homogeneous of degree $|V|$ with spectrum $W/V$ of dimension $n$; it has nuclear dimension $n$, and nuclear dimension does not increase on passing to ideals \cite{WZ2010}. So $\dimnuc(\Cr(D_\infty^n))=n$, $\dimnuc(\Cr(\Gamma_0))=3$ and $\dimnuc(\Cr(\Gamma_4))=4$. The group $\mathrm{SL}(3,\mathbb Z)$ is exact, being linear \cite{GHW2005}, and its reduced group algebra does not have the completely bounded approximation property \cite{H2016}. Thus $\rdim$ does not detect weak amenability.

\begin{remark}[Relative dimension versus coarse dimension]
The original comparison with asymptotic dimension is strictly different for the relative invariant. For every $n\ge1$, the group $D_\infty^n$ has asymptotic dimension $n$, its uniform Roe algebra has nuclear dimension one, and its regular inclusion has nuclear dimension zero. In particular, even within virtually abelian groups the regular inclusion need not recover the dimension of the uniform Roe algebra. Conversely, $\mathbb Z^n$ has relative nuclear dimension one, no matter how large $n$ is.
\end{remark}

\begin{remark}[The weak$^*$ variant and a different ambient algebra]
Proposition~\ref{prop:obstruction} shows that a finite-colour weak$^*$ approximation of the canonical inclusion $A\to A^{**}$, with the uniform bounds specified there, forces nuclearity of $A$. We do not assert the converse for a dimension-zero weak$^*$ variant without specifying its approximation requirements. On the other hand, if $A$ is a separable exact C$^*$-algebra, the embedding theorem of Kirchberg and Phillips \cite{KP2000} gives an inclusion into $\mathcal O_2$, and $\dimnuc(\mathcal O_2)=1$. Hence the infimum of the nuclear dimension of inclusions of $A$ into nuclear C$^*$-algebras is at most one. This is a different question from the canonical ambient algebra used here.
\end{remark}

\begin{question}\label{q:classification}
Which countable amenable groups with torsion abelianization have relative nuclear dimension zero? In particular, what is the value for infinite finitely generated amenable torsion groups?
\end{question}
The torsion assumption on the abelianization is necessary by Theorem~\ref{thm:characters}, and it is not sufficient by Theorem~\ref{thm:holonomy} and Example~\ref{ex:orientation}. Theorem~\ref{thm:dihedral} shows that the absence of infinite-order elements in the group itself is not necessary. The remaining case includes, for example, the first Grigorchuk group; no value is claimed here for such groups.

\begin{question}\label{q:kcriterion}
For an amenable group $\Gamma$, is $K_1(\iota_\Gamma)=0$ sufficient for $\rdim(\Gamma)=0$? If not, what additional obstructions arise from ordered $K$-theory or from other invariants of zero-dimensional morphisms?
\end{question}
The implication from dimension zero to vanishing of $K_1(\iota_\Gamma)$ is established by Proposition~\ref{prop:kone}. Proposition~\ref{prop:kone} also gives the vanishing of $K_1(\iota_\Gamma;\mathbb Z/q\mathbb Z)$ for $q\ge2$, and Castillejos and Neagu obtain further obstructions in terms of the total invariant \cite{CN2024}. For the groups $D_\infty^n$ the coefficient condition holds for a simple reason: by the interval-algebra description of $\Cr(D_\infty)$ given above, there is an exact sequence $0\to C_0((0,1),M_2)\to\Cr(D_\infty)\to\mathbb C^4\to0$, the quotient map recording the diagonal values at the two endpoints. The index map $K_0(\mathbb C^4)=\mathbb Z^4\to K_1(C_0((0,1),M_2))=\mathbb Z$ is, up to sign, $(a,b,c,d)\mapsto a+b-c-d$, which is surjective, and the six-term exact sequence gives $K_0(\Cr(D_\infty))\cong\mathbb Z^3$ and $K_1(\Cr(D_\infty))=0$. As $\Cr(D_\infty)$ is type I, the K\"unneth theorem gives $K_0(\Cr(D_\infty^n))\cong\mathbb Z^{3^n}$, which is torsion-free, and $K_1(\Cr(D_\infty^n))=0$; the exact sequence in the proof of Proposition~\ref{prop:kone} then gives $K_1(\Cr(D_\infty^n);\mathbb Z/q\mathbb Z)=0$. No sufficiency statement is made here.

\begin{question}\label{q:virtuallyabelian}
Let $\Gamma$ be a finitely generated group with a normal subgroup $H\cong\mathbb Z^n$ of finite index and holonomy group $F$. If $(\Lambda^{2k+1}\mathbb Q^n)^F=0$ for all $k\ge0$, is $\rdim(\Gamma)=0$?
\end{question}
By Theorem~\ref{thm:holonomy}, the condition in Question~\ref{q:virtuallyabelian} is equivalent to $K_1(\iota_\Gamma)=0$, and it is necessary for $\rdim(\Gamma)=0$ by Proposition~\ref{prop:kone}. So Question~\ref{q:virtuallyabelian} is the virtually abelian case of Question~\ref{q:kcriterion}. The groups $D_\infty^n$ satisfy the condition and have dimension zero, by Theorem~\ref{thm:dihedral} and Corollary~\ref{cor:products}. The permanence results of Proposition~\ref{prop:permanence} give further groups of dimension zero, such as $D_\infty^n\times L$ with $L$ locally finite; when $L$ is infinite these groups are not finitely generated, so they lie outside the scope of the finite-rank holonomy criterion. By Corollary~\ref{cor:bieberbach}, every group satisfying the condition has torsion. A natural approach to the general case is to look for actions on the Cantor set whose reduced crossed products are AF, as in Theorem~\ref{thm:dihedral}.

\begin{question}\label{q:AF}
If $\rdim(\Gamma)=0$, does there exist a separable unital AF algebra $B$ with $\Cr(\Gamma)\subset B\subset\Cu(\Gamma)$? Can such a factorization be realized by a compact metrizable $\Gamma$-space with an AF reduced crossed product?
\end{question}
Lemma~\ref{lem:localfd} only supplies local finite-dimensional approximations and does not assert that they can be chosen as a nested family. Castillejos and Neagu characterize zero-dimensional $*$-homomorphisms, within the framework of the classification of $*$-homomorphisms, as those that approximately factor through an AF algebra \cite{CN2024}. Question~\ref{q:AF} asks for more: an actual AF subalgebra lying between $\Cr(\Gamma)$ and $\Cu(\Gamma)$. The AF intermediate algebra in Theorem~\ref{thm:dihedral} is supplied by a specific dynamical construction.

\begin{question}\label{q:propagation}
Can one refine $\rdim$ by controlling the propagation of the finite-dimensional approximations as a function of the tolerance and the finite set of group elements?
\end{question}
For example, in dimension zero one can ask for bounds on the propagation needed to approximate the unit ball of a finite-dimensional algebra $F$ as in Lemma~\ref{lem:localfd}. Since $F$ is finite-dimensional, its unit ball is compact and finite-propagation operators are dense in $\Cu(\Gamma)$, so a common finite propagation bound exists for each fixed tolerance. A quantitative invariant would need to control its dependence on the generating set and tolerance. The twisted partitions of unity in \cite{GLZZ2026} supply a natural two-colour comparison. However, Proposition~\ref{prop:Zindex} precludes a one-colour construction for $\mathbb Z$ in the canonical uniform Roe algebra, even though the regular representation is quasidiagonal by amenability \cite{TWW2017}. The dihedral example shows that allowing a change of the group and its regular inclusion can remove this obstruction.


\end{document}